\documentclass[a4paper,12pt]{article}
\usepackage{amsmath}
\usepackage{amsfonts}
\usepackage{amssymb}
\usepackage{latexsym}
\usepackage{epsfig}
\usepackage{graphicx}
\usepackage{oldgerm}
\usepackage{theorem}
\def\C{\mathbb{C}}
\def\D{\mathbb{D}}
\def\E{\mathbb{E}}
\def\I{\mathbb{I}}
\def\N{\mathbb{N}}
\def\P{\mathbb{P}}
\def\R{\mathbb{R}}
\def\X{\mathbf{X}}
\def\Y{\mathbf{Y}}
\def\Z{\mathbb{Z}}
\def\f{\mathbf{f}}
\def\t{\mathbf{t}}
\def\x{\mathbf{x}}
\def\1{\mathbf{1}}
\def\mM{\mathfrak{M}}
\def\mX{\mathfrak{X}}
\def\u{\mathfrak{u}}
\def\cA{\mathcal{A}}
\def\cD{\mathcal{D}}
\def\cE{\mathcal{E}}
\def\cP{\mathcal{P}}
\def\Ai{\mathrm{Ai}}

\def\Det{\mathrm{Det}}
\def\supp{\mathrm{supp}\ }
\theorembodyfont{\itshape}
\newtheorem{thm}{Theorem}[section]
\newtheorem{lem}[thm]{Lemma}
\newtheorem{cor}[thm]{Corollary}
\newtheorem{prop}[thm]{Proposition}

\newcommand{\SSC}[1]{\section{#1}\setcounter{equation}{0}}
\newcommand{\qed}{\hbox{\rule[-2pt]{3pt}{6pt}}}

\begin{document}
\title{Strong Markov property of determinantal processes \\
with extended kernels}

\author{
Hirofumi Osada
\footnote{
Faculty of Mathematics, Kyushu University,
Fukuoka 819-0395 Japan;
e-mail: osada@math.kyushu.ac.jp
}
\quad 
Hideki Tanemura
\footnote{
Department of Mathematics and Informatics,
Faculty of Science, Chiba University, 
1-33 Yayoi-cho, Inage-ku, Chiba 263-8522, Japan;
e-mail: tanemura@math.s.chiba-u.ac.jp
}}

\pagestyle{plain}
\maketitle
\begin{abstract}
Noncolliding Brownian motion (Dyson's Brownian motion model with parameter $\beta=2$)
and noncolliding Bessel processes are determinantal processes; that is,
their space-time correlation functions are represented by determinants. Under a proper scaling limit,
such as the bulk, soft-edge and hard-edge scaling limits, these processes converge to determinantal processes describing systems with an infinite number of particles. The main purpose of this paper is to show the strong Markov property of these limit processes, which are determinantal processes with the extended sine kernel, extended Airy kernel and extended Bessel kernel, respectively. We also determine the quasi-regular Dirichlet forms and infinite-dimensional stochastic differential equations associated with the determinantal processes. 



\end{abstract}

\SSC{Introduction}

In a system of $N$ independent one-dimensional diffusion processes, if we impose the condition that the particles never collide with one another, we obtain an interacting particle system with a long-range strong repulsive force between any pair of particles. We call such a system a noncolliding diffusion process.  
In 1962 Dyson \cite{Dys62} showed that when the individual diffusion process is one-dimensional Brownian motion (BM), it is related to a matrix-valued process. We call this stochastic process {\it the Dyson model}. The model solves the stochastic differential equation (SDE)
\begin{equation}
dX_j(t)=dB_j(t)+ \sum_{\substack{k:\, k \not= j \\ k = 1 }}^{N} 
\frac{dt}{X_j(t)-X_k(t)},
\quad j \in\I_N\equiv \{1,2,\dots,N\},
\label{eqn:Dyson1}
\end{equation}
where $B_j(t), j=1,2,\dots,N$ are independent one-dimensional BMs.
When the individual diffusion process is a squared Bessel process with index $\nu>-1$, 
the noncolliding diffusion process is called a {\it noncolliding squared Bessel process}, satisfying the SDE
\begin{equation}\label{eqn:Bessel1}
dZ_j^{\nu}(t)=2\sqrt{Z_j^{\nu}(t)}dB_j(t)+2(\nu+1)dt
+ \sum_{\substack{k:\, k \not= j \\ k = 1 }}^{N} \frac{4 Z_j^{\nu}(t) \ ds}{Z_j^{\nu}(t)-Z_k^{\nu}(t)},
\ j\in \I_N
\end{equation}
and if $-1 < \nu < 0$ having a reflection wall at the origin.
These processes dynamically simulate the eigenvalue statistics
of the Gaussian random matrix ensembles studied in
random matrix theory \cite{Meh04,For10}.

Let $\mM$ be the space of nonnegative integer-valued Radon measures on $\R$.
This space is a Polish space with the {\it vague topology}.
The space can be regarded as a configuration space of unlabeled particles in $\R$.
A probability measure on the space $\mM$ is called a {\it determinantal point process} (DPP) or {\it Fermion point process},
if its correlation functions are generally represented by determinants \cite{ST03, Sos00}.
In this paper we say that an $\mM$-valued process $\Xi(t)$ is {\it determinantal} 
if the multitime correlation functions for any chosen series of times are represented by determinants. It has been shown that for any initial configuration $\xi^N= \sum_{j=1}^N \delta_{x_j}$,
the unlabeled process $\Xi^N(t) = \sum_{j=1}^N \delta_{ X_j(t)}$ is a determinantal process for the Dyson model in \cite{KT10a}, and for the noncolliding Bessel process in \cite{KT11a}.

Suppose that $(\Xi(t),\P^{\xi^N})$ is a noncolliding diffusion process starting from the initial configuration $\xi^N$  of a finite number of particles. Let $\xi$ be the configuration of an infinite number of particles, and $\xi_{[-L,L]}$ denote the the restriction of $\xi$ on the set $[-L,L]$.
When the sequence of processes $(\Xi(t), \P^{\xi_{[-L,L]}})$, $L\in\N=\{1,2,\dots\}$, converges to an $\mM$-valued process starting from the configuration $\xi$,
we denote the limit process by $(\Xi(t), \P^\xi)$,
and say that process $(\Xi(t), \P^\xi)$ is well-defined.
Sufficient conditions were given for configuration $\xi$ that limit process is well definied, for noncolliding BM with bulk scaling \cite{KT10a},  noncolliding BM with soft edge scaling \cite{KT09}, and noncolliding Bessel processes with hard-edge scaling \cite{KT11a}.
We denote the assocoated limit processes by $(\Xi(t), \P_{\sin}^\xi)$, $(\Xi(t), \P_{\Ai}^\xi)$ and $(\Xi(t), \P_{\nu}^\xi)$, respectively.

Let $\mu_{\sin}$, $\mu_{\Ai}$ and $\mu_{J_\nu}$ be
the DPPs with the sine kernel (\ref{S-kernel}), the Airy kernel (\ref{A-kernel}) and 
the Bessel kernel (\ref{B-kernel}), respectively.
They are the probability measures
obtained in the bulk scaling limit and soft-edge scaling limit
of the eigenvalue distribution (\ref{GUE_N}) in the 
{\it Gaussian unitary ensemble} (GUE),
and in the hard-edge scaling limit of the eigenvalue distribution (\ref{chGUE_N}) in the
{\it chiral Gaussian unitary ensemble} (chGUE),
respectively.
It was shown in \cite{KT11b} that processes $(\Xi(t), \P_{\sin}^\xi)$, $(\Xi(t), \P_{\Ai}^\xi)$ and  $(\Xi(t), \P_{\nu}^\xi)$ have DPPs $\mu_{\sin}$, $\mu_{\Ai}$ and $\mu_{J_\nu}$ as reversible measures, and the reversible processes coincide with determinantal processes $(\Xi(t), \P_{\sin})$ with the extended sine kernel (\ref{ex_S-kernel}), $(\Xi(t), \P_{\Ai})$ with the extended Airy kernel (\ref{ex_A-kernel}), and $(\Xi(t), \P_{\nu})$ with
 the extended Bessel kernel (\ref{ex_B-kernel}), respectively.
The main purpose of this paper is to prove the
{\it strong Markov property} of these infinite-dimensional
determinantal processes $(\Xi(t) ,\P_{\star})$, $\star\in\{\sin, \Ai\} \cup \{J_\nu\ ;\nu>-1\}$ (Theorem \ref{TH:SM}).

For each $\star\in\{\sin, \Ai\} \cup \{J_\nu\ ;\nu>-1\}$,
the diffusion process $(\Xi(t) ,\mathbf{P}_{\star}^\xi)$  associated with DPP $\mu_{\star}$ was constructed by a Dirichlet form technique presented by the first author \cite{o.rm, o.rm2}, and associated infinite-dimensional stochastic differential equation (ISDE) was derived \cite{o.isde, ot.1, ho.1}. The relation between processes $(\Xi(t), \P_\star ^\xi)$ and $(\Xi(t), \mathbf{P}_\star ^\xi)$ was first discussed in \cite{KT07a}; it was later proven that both of these are extensions of the same pre-Dirichlet form on the set of polynomial functions \cite{ot.0}. However, the coincidence of these two processes has long been an open problem. Using Theorem \ref{TH:SM}, it can be proved affirmatively (Theorem \ref{TH:SDE}). 
The coincidence of processes from completely different approaches of construction enable us to examine them from various points of view.
From the algebraic construction by virtue of their determinantal structure of their space-time correlation functions, we obtain  {\it quantitative} information of the processes such as the moment generating functions of multitime distributions (\ref{eqn:Fred}) given by the Fredholm determinant of space-time correlation functions \cite{KT09, KT10a, KT11a}; while from the analytic construction through Dirichlet form theory, we deduce many {\it qualitative} properties of sample paths of the labeled diffusion by means of the ISDE representation, which makes us possible to apply the Ito calculus to the processes \cite{Osa04, o.isde, ot.1}.

The proof is based on the following two facts:

\vskip 1mm

\noindent (i) Determinantal process $(\Xi(t), \P_\star^\xi)$
solves the same ISDE as process $(\Xi(t), \mathbf{P}_\star^\xi)$.

\noindent (ii) The solutions of the ISDE are  unique.

\vskip 1mm

\noindent Fact (ii) was shown in \cite{ot.2}.
Fact (i) is derived by the arguments in \cite{o.isde}, in which the consistent families of diffusion processes $(\X^k(t), H(t))$, $k\in \N$ describing the joint processes of tagged particles and unlabeled infinite particles are introduced. To construct the processes, the strong Markov property of the unlabeled process $(\Xi(t), \P_\star^\xi)$ is crucial \cite{o.tp}.

To prove the strong Markov property of the determinantal process $(\Xi(t), \P_\star^\xi)$ 
we introduce the subsets $ \mX_{\star}$ in \eqref{:21}--\eqref{:23},
and equip them with the topology defined by the inductive limit. 
One of the key points for the strong Markov property is to prove the Feller-like property of the process 
$(\Xi(t), \P_\star ^\xi)$ restricted on $ \mX_{\star}$.
For showing the Feller property of the process it is necessary to have if $\xi_n\to\xi$ in $\mX_{\star}$ with the topology, $\P^{\xi_n} \to \P^{\xi}$ weakly on $C([0,\infty),\mX_\star)$.
In stead of the property, 
we prove that  $\P_{\star}^\xi (\Xi \in C([0,\infty],\mX_\star))=1$ for $\mu_{\star}$-a.s. $\xi$ in Proposition \ref{Prop:3.7},
and combine it with the fact that if $\xi_n\to\xi$ in $\mX_{\star}$ with the topology, $\P^{\xi_n} \to \P^{\xi}$ weakly on $C([0,\infty],\mX)$, which was already proved in \cite{KT09, KT10a, KT11a}, and derive the strong Markov property of the revrsible process $(\Xi(t), \P_\star)$.

Let $\mu_{\sin}^N$, $\mu_{\Ai}^N$ and $\mu_{J_\nu}^N$ be probability measures on the configuration spaces of finite numbers of particles defined by (\ref{mu_sin_N}), (\ref{mu_Ai_N}) and (\ref{mu_J_N}), respectively.
 For each $\star\in\{\sin, \Ai\} \cup \{J_\nu\ ;\nu>-1\}$, $\mu_\star$ is obtained by the limit of $\mu_{\star}^N$, $N\to \infty$.
It is interesting to study the relation between the Dirichlet forms related to $\mu_{\star}^N$ and $\mu$ with the same square field defined in (\ref{:sqfield}).  We denote by $(\Xi(t), \mathbf{P}_{\star}^{\xi^N})$ the Diffusion process associated the Dirichlet form related to $\mu_{\star}^N$ and starting from any configuration of $N$ particles.
 As one example of byproducts of Theorem \ref{TH:SDE}, we obtain that $(\Xi(t), \mathbf{P}_{\star}^{\xi_N^{\sf l}})$ converges to $(\Xi(t), \mathbf{P}_{\star}^\xi$), $N\to\infty$, where $\xi_N^{\sl l}$ is the configuration given in (\ref{:xil}) (Corollary \ref{C:23}).

The paper is organized as follows. In Section 2 preliminaries and main results are given. In section 3 we give the proof of the main theorems.

\vskip 10mm

\SSC{Preliminaries and main results}

Let $\mM = \mM (\R)$ be the space of nonnegative 
integer-valued Radon measures on $\R$,
which is a Polish space with the vague topology.
We say $\xi_n, n\in\N$ converges to $\xi$ vaguely if 
$\lim_{n \to \infty} \int_{\R} \varphi(x) \xi_n(dx)
=\int_{\R} \varphi(x) \xi(dx)$ 
for any $\varphi \in C_0(\R)$, where $C_0(\R)$ is the set of all 
continuous real-valued functions with compact supports. A subset ${\mathfrak N}$ of $\mM$ is relatively compact if and only if $\sup_{\xi\in\mathfrak N}\xi (K) <\infty$, for any compact set $K$ in $\R$. Each element $\xi$ of $\mM$ can be represented as
$\xi(\cdot) = \sum_{j\in \Lambda}\delta_{x_j}(\cdot)$
with an index set $\Lambda$ and a sequence of points 
$\x =(x_j)_{j \in \Lambda}$ in $\R$ 
satisfying $\xi(K)=\sharp\{x_j: x_j \in K\} < \infty$ 
for any compact subset $K \subset \R$.
We introduce the subspace $\mM_0$ of $\mM$ defined as $\mM_0= \{ \xi\in\mM: \xi(\{x\})\le 1, \quad x\in\R \}$. Any element $\xi\in\mM_0$ is identical to its support,
and can be regarded as a countable subset $\{x_j\}_{j\in\Lambda}$.
For an element $\xi\in \mM\setminus \mM_0$, there is a point $x\in\supp \xi$
such that $\xi(\{x\})\ge 2$. Such a point is called a multiple point.
We call an element $\xi$ of $\mM$ an {\it unlabeled} configuration, and a sequence of points $\x=(x_j)_{j\in\Lambda}$ a {\it labeled} configuration.
For $\xi(\cdot)=\sum_{j\in \Lambda}
\delta_{x_j}(\cdot) \in\mM$,
we introduce the following operations
\begin{description}

\item[(restriction)] 
$\xi_A =\xi \cap A =\displaystyle{\sum_{x \in A}} 
\delta_{x}(\cdot)$, \ for $A\subset \R$,

\item[(shift)]  
$\tau_u \xi(\cdot) =\displaystyle{\sum_{j \in \Lambda}} 
\delta_{x_j+u}(\cdot)$, \ for $u \in \R$. 
\end{description}
We also introduce the configuration space ${\mM}^+ = \{\xi\in \mM: \xi((-\infty,0))=0\}$.

A probability measure on the space $\mM$ is called a DPP or {\it Fermion point process},
if its correlation functions are generally represented by determinants \cite{ST03, Sos00}. In this paper we consider the DPPs with the following three continuous kernels on $\R\times\R$.

\vskip 3mm

\noindent (i) Sine kernel $K_{\sin}(x,y)$:
\begin{equation}\label{S-kernel}
K_{\sin}(x,y)=
\frac{1}{2 \pi} \int_{|k| \leq 1} dk \,
e^{\sqrt{-1} k(x-y)}
= \left\{ \begin{array}{ll}
\displaystyle{\frac{\sin (x-y)}{\pi (x-y)}}
&\mbox{if } x \not= y, \cr
\cr
\displaystyle{\frac{1}{\pi}}
&\mbox{if } x=y. \end{array} \right.
\end{equation}

\noindent (ii) Airy kernel $K_{\Ai}(x,y)$:
\begin{equation}\label{A-kernel}
K_{\Ai}(x,y) = \left\{ \begin{array}{ll}
\displaystyle{
\frac{\Ai(x) \Ai'(y)-\Ai'(x) \Ai(y)}{x-y}}
&\mbox{if } x \not= y, \cr
\cr
(\Ai'(x))^2-x (\Ai(x))^2
&\mbox{if } x=y, \end{array} \right.
\end{equation}
where $\Ai(x) = \frac{1}{2 \pi} \int_{\R} dk e^{\sqrt{-1}(x k+k^3/3)}$, $x\in\R$, is the Airy function
and $\Ai'(x)=d \Ai(x)/dx$.

\noindent (iii) Bessel kernel $K_{J_\nu}(x,y)$:
\begin{equation}\label{B-kernel}
K_{J_\nu}(x,y) = \left\{ \begin{array}{ll}
\displaystyle{
\frac{J_{\nu}(2 \sqrt{x}) \sqrt{y} J_{\nu}'(2 \sqrt{y})
-\sqrt{x} J_{\nu}'(2\sqrt{x}) J_{\nu}(2\sqrt{y})}{x-y}} 
& \mbox{if} \quad x \not= y, \cr
\cr
 J_{\nu}(2 \sqrt{x})^2-J_{\nu+1}(2 \sqrt{x}) J_{\nu-1}(2 \sqrt{x})
& \mbox{if} \quad x=y,
\end{array} \right.
\end{equation}
where $J_{\nu}(\cdot)$ is the Bessel function with index $\nu>-1$ defined as
\begin{equation}
J_{\nu}(z) = \sum_{n=0}^{\infty} 
\frac{(-1)^n}{\Gamma(n+1) \Gamma(n+1+\nu)}
\left( \frac{z}{2} \right)^{2n+\nu},
\quad z\in \C
\nonumber
\end{equation}
with the gamma function 
$\Gamma(x)=\int_0^{\infty} e^{-u} u^{x-1} du$ and $J_{\nu}'(x)=dJ_{\nu}(x)/dx$.
We denote by $\mu_{\sin}$, $\mu_{\Ai}$, $\mu_{J_\nu}$, DPPs with the sine, Airy and Bessel kernels, respectively.
We note that DPPs $\mu_{\sin}$, $\mu_{\Ai}$ and $\mu_{j_\nu}$ are the limit distributions of DPPs $\mu_{\sin}^N$, $\mu_{\Ai}^N$, and $\mu_{J_\nu}^N$ 
defined, respectively, in (\ref{mu_sin_N}), (\ref{mu_Ai_N}), (\ref{mu_J_N}), which are related to the eigenvalue distributions of random matrices \cite{Meh04, AGZ10}.

Let $\rho(dx)=\rho(x)dx$ be a Radon measure on $\R$ with density $\rho(x)$ such that $\rho(\R)=\int_\R \rho(x)dx\in \{0\}\cup \N \cup \{\infty\}$. For $\varepsilon, \kappa >0$ and $m_0, L_0\in\N$, we denote by $\mX_{L_0, m_0}^{\varepsilon, \kappa}(\rho)$ the set of configurations $\xi\in \mM$ satisfying $\xi(\R)=\rho(\R)$,
\begin{equation}\nonumber
\left|\rho([0,L]) -\xi([0,L])\right|\le L^{\varepsilon}, \;
\left|\rho([-L,0]) -\xi([-L,0])\right|\le L^{\varepsilon},
\quad L\ge L_0,
\end{equation}
and 
\begin{equation}\nonumber
\max_{k\in\Z}\xi\bigg( [g^\kappa(k), g^\kappa(k+1)]\bigg)\le m_0,
\end{equation}
where 
\begin{equation}\label{gk}
g^\kappa (x)={\rm sign}(x)|x|^\kappa, \ x\in\R.
\end{equation}
The set $\mX_{L_0, m_0}^{\varepsilon, \kappa}(\rho)$ is then relatively compact with the vague topology.
We introduce the following three configurations associated with DPPs $\mu_\star$, $\star \in \{\sin, \Ai\} \cup \{J_\nu\ ;\nu>-1\}$:
\begin{eqnarray} \label{:21}
&&\mX_{\sin} =\bigcup_{L_0, m_0\in\N}\bigcup_{\substack{\varepsilon\in (0,1)\\ \kappa\in (0,\kappa_{\sin})}}\mX_{L_0,m_0}^{\varepsilon,\kappa}(\rho_{\sin}),
\\\label{:22}
&&\mX_{\Ai} =\bigcup_{L_0, m_0\in\N}\bigcup_{\substack{\varepsilon\in (0,1)\\ \kappa\in (0,\kappa_{\Ai})}}\mX_{L_0,m_0}^{\varepsilon, \kappa}(\rho_{\Ai}),
\\\label{:23}
&&\mX_{J_\nu} =\mM^+ \cap\bigcup_{L_0, m_0\in\N}\bigcup_{\substack{\varepsilon\in (0,1)\\ \kappa\in (0,\kappa_{J_\nu})}}\mX_{L_0,m_0}^{\varepsilon, \kappa}(\rho_{J_\nu}),
\end{eqnarray}
where $\kappa_{\sin}=1$, $\kappa_{\Ai}=2/3$, $\kappa_{J_\nu}=2$, and $\rho_{\star}$ is the first order correlation function of DPPs $\mu_{\star}$ for each $\star \in \{\sin, \Ai\} \cup \{J_\nu\ ;\nu>-1\}$. 
It is readily verified that $\mu_{\star}(\mX_{\star})=1$, for $\star =\sin$, $\Ai$ or $J_\nu$. We equip space $\mX_{\star}$ with the topology defined by the inductive limit. Note that this is a locally compact space with the topology. We denote by $C([0,\infty),\mX_\star)$ the set of $\mX_\star$-valued continuous functions with the topology. Then $\Xi(\cdot)\in C([0,\infty),\mX_\star)$ implies that $\Xi(\cdot)$ is continuous with the vague topology, and
for any $T>0$ there exist $M_0,L_0 \in\N$, $\varepsilon \in (0,1)$ and  $\kappa \in (0,\kappa_*)$ such that $\Xi(t)\in \mX_{L_0,m_0}^{\varepsilon,\kappa}(\rho_*)$ for $t\in [0,T]$. 

In this paper we say that an $\mM$-valued process $\Xi(t)$ is {\it determinantal}, if the multitime correlation functions for any chosen series of times are represented by determinants. 
In other words, a {\it determinantal process} is 
an $\mM$-valued process such that, for any integer $M\in \N$, $\f=(f_1,f_2,\dots,f_M) \in C_0(\R)^M$, and sequence of times $\t=(t_1,t_2,\dots,t_M)$ with $0 < t_1 < \cdots < t_M < \infty$, the moment generating function of multitime distribution, 
$
{\Psi}^{\t}[\f]
\equiv \E \left[\exp \left\{ \sum_{m=1}^{M} 
\int_{\R} f_m(x) \Xi(t_m, dx) \right\} \right],
$
is given by a Fredholm determinant
\begin{equation}
{\Psi}^{\t}[\f]
= \mathop{\Det}_
{\substack
{(s,t)\in \{t_1, t_2, \dots, t_M \}^2 \\
(x,y)\in \R^2}}
\Big[\delta_{s t} \delta(x-y)
+ \mathbb{K}(s, x; t, y) \chi_{t}(y) \Big],
 \label{eqn:Fred}
\end{equation}
where $\chi_{t_m}= e^{f_m}-1$, $1\le m\le M$, and $\mathbb{K}$ is a locally integrable function  called a (space-time) {\it correlation kernel} \cite{KT07a, KT09, KT10a, KT11a}.
In this paper we study the determinantal processes with the following three correlation kernels:
\vskip 3mm

\noindent (i){\it Extended sine kernel} 
${\mathbb{K}}_{\sin}(s,x;t,y), s,t\in\R^+\equiv \{x \in \R : x \geq 0 \}, x,y\in\R$:
\begin{equation}\label{ex_S-kernel}
{\mathbb{K}}_{\sin}(s,x;t,y) 
= \left\{ \begin{array}{ll} 
\displaystyle{\frac{1}{\pi}
\int_{0}^{1} du \, e^{ u^2 (t-s)/2} 
\cos \{ u (y-x)\} }
& \mbox{if $s<t, $} \cr
K_{\sin}(x,y)
& \mbox{if $s=t$,} \cr
\displaystyle{
- \frac{1}{\pi}\int_{1}^{\infty} du \, 
e^{ u^2 (t-s)/2} \cos \{ u (y-x) \} }
& \mbox{if $s>t$}.
\end{array} \right.
\end{equation}

\noindent (ii) {\it Extended Airy kernel} ${\mathbb{K}}_{\Ai}(s,x;t,y), s,t\in\R^+, x,y\in\R$:
\begin{equation}\label{ex_A-kernel}
{\mathbb{K}}_{\Ai}(s,x;t,y) = \left\{ 
\begin{array}{ll}
\displaystyle{
\int_{0}^{\infty} d u \, e^{-u(t-s)/2} \Ai(u+x) \Ai(u+y)}
& \mbox{if $s < t,$} \cr
K_{\Ai}(x,y) & \mbox{if $s =t,$}\cr
\displaystyle{- \int_{-\infty}^{0} d u \, e^{-u(t-s)/2} \Ai(u+x) \Ai(u+y)}
& \mbox{if $s>t$}.
\end{array}
\right.
\end{equation}

\noindent(iii) 
{\it Extended Bessel kernel} $\mathbb{K}_{J_{\nu}}(s,x;t,y), s,t\in\R^+, x,y\in\R^+$:
\begin{equation}\label{ex_B-kernel}
\mathbb{K}_{J_{\nu}}(s,x;t,y) = \left\{
   \begin{array}{ll}
\displaystyle{
\int_{0}^{1} d u \,
e^{-2u(s-t)} J_{\nu}(2 \sqrt{u x})
J_{\nu}(2 \sqrt{u y})
} 
& \mbox{if} \quad s < t,  \\
\displaystyle{
K_{J_\nu}(x,y)}& \mbox{if} \quad s=t, \\
\displaystyle{
- \int_{1}^{\infty} d u \,
e^{-2u(s-t)} J_{\nu}(2 \sqrt{ux})
J_{\nu}(2 \sqrt{uy})
}
& \mbox{if} \quad s > t.
   \end{array} \right. 
\end{equation}

We denote the determinantal process with the extended sine kernel $\mathbb{K}_{\star}(s,x;t,y)$ by $(\Xi(t), \P_{\star})$, $\star=\sin, \Ai, J_\nu$. 
The determinantal process $(\Xi(t), \P_{\star})$ is reversible with the DPP $\mu_{\star}$, for $\star=\sin, \Ai$ and $J_\nu$. 

The main theorem of this present paper is as follows.
\begin{thm}
\label{TH:SM}
Let $\star \in \{\sin, \Ai\} \cup \{J_\nu\ ;\nu>-1\}$.
The determinantal process $(\Xi(t), {\P}_{\star})$
is a reversible diffusion process.
\end{thm}

Diffusion processes with reversible probability measures $\mu_\star, \star\in\{\sin, \Ai\} \cup \{J_\nu\ ;\nu>-1\}$ were constructed by means of the Dirichlet form technique by the first author \cite{o.rm, o.rm2} and associated ISDEs were derived in \cite{o.isde, ot.1, ho.1}. The second theorem of this paper is the coincidence of the process constructed by the Dirichlet form technique and the determinantal process $(\Xi(t), {\P}_{\star})$, $\star \in \{\sin, \Ai\} \cup \{J_\nu\ ;\nu>-1\}$.
To state the theorem we prepare some notation and review the related results.
A function $f$ on $\mM$ is called {\it local} if $f(\xi)=f(\xi_K)$ for some compact set $K$.  
For a local function $f$ with $f(\xi)=f(\xi_K)$ 
we introduce the functions $\check{f}_k$ 
on $\R^k$, $k\in\N_0\equiv \{0\}\cup \N $ defined as
$\check{f}_0=f(\emptyset)$, where $\emptyset$ is the null configuration, and for $ k\in \N $, as
$$
\check{f}_k(\x_k)=f\left(\sum_{j=1}^k \delta_{x_j} \right)   
\text{ for }
\x_k \in K^k 
.$$
We extend the domain of $ \check{f}_k(\x_k)$ to $ \R^k\backslash K^k$ 
by the consistency arising from $ f (\xi) = f (\xi_K)$. 
Hence $\check{f}_k$, $k\in\N_0$, satisfy the consistency relation
\begin{equation}
\nonumber 
\check{f}_{k+1}(\x_k, y)=\check{f}_{k}(\x_k), \ \x_k\in \R^k, \ y\notin K.
\end{equation}
A local function $f$ is called {\it smooth} if the $\check{f}_k$ are smooth 
for $k\in\N_0$. 
We denote by ${\cal D}_{\infty}$ the set of all local smooth functions on $\mM$.

Let $\x_k=(x_i)_{i=1}^k \in \R^k$, $k\in\N_0$, and
$f,g \in {\cal D}_\infty$; then
\begin{equation}\label{:sqfield}
{\Bbb D}^a(f,g)(\x_k)= \frac{1}{2}
\sum_{i=1}^k a(x_i)
\frac{\partial \check{f}_k(\x_k)}{\partial x_{i}}
\frac{\partial \check{g}_k(\x_k)}{\partial x_{i}}
,
\end{equation}
where $a$ is a smooth function on $\R$. 
For given $f,g\in{\cal D}_\infty$, the right-hand side is a permutation invariant function, and the square field ${\Bbb D}(f,g)$ can be regarded as a local function with variable $\xi=\sum_{i\in \N}\delta_{x_i}\in \mM$.
For probability $\mu$ on $\mM$,
we denote by $L^2(\mM,\mu)$ the space of square integrable functions on $\mM$ with
the inner product $\langle \cdot, \cdot \rangle_\mu$ and the norm $\|\cdot\|_{L^2(\mM,\mu)}$.
We consider the bilinear form $(\cE_\mu^a, \cD_\infty^\mu)$ on $L^2(\mM,\mu)$ defined as
\begin{eqnarray}
&&\cE_\mu^a (f,g)= \int_{\mM}{\Bbb D}^a(f,g)d\mu,
\nonumber 
\\
&&\cD_\infty^\mu=\{f\in \cD_\infty : \|f\|_1^2 <\infty \}
\nonumber 
,\end{eqnarray}
where $\|f\|_1^2 \equiv \cE_\mu^a (f,f)+\|f \|_{L^2(\mM,\mu)}^2$. 
We consider three cases $(\mu, a(x))= (\mu_{\sin}, 1), (\mu_{\Ai},1)$ and  $(\mu_{J_\nu}, 4x)$.
For each case it is proved in \cite{o.rm, o.rm2} that the closure $(\cE_{\star}, \cD_{\star})$ of the bilinear form  $(\cE_{\mu_{\star}}^a, \cD_\infty^{\mu_\star})$ is a quasi-Dirichlet form, associated with the diffusion process
$(\Xi(t), \mathbf{P}_\star^\xi)$, $\star\in\{\sin, \Ai\} \cup \{J_\nu\ ;\nu>-1\}$.
One can naturally lift each unlabeled path $\Xi$ to the labeled path $\X$ by a labeled map $\mathfrak{l}=(\mathfrak{l})_{i\in\N}$ (see for example \cite[Theorem 2.4]{o.tp}).
It was also proved that a labeled process $\X =(X_j)_{j\in\N}$ with $\Xi(t)= \sum_{j\in\N}\delta_{X_j(t)}$ 
solves the ISDE 
\begin{equation}\label{ISDE:sin}
\tag{sin}
dX_j(t)=dB_j(t)+\lim_{r\to\infty}\sum_{\substack{k:\, k\not=j \\ |X_k(t)|<r}}^\infty\frac{dt}{X_j(t)-X_k(t)},
\end{equation}
under $\mathbf{P}_{\sin}^\xi$, 
  a labeled process $\Y =(Y_j)_{j\in\N}$ with 
$\Xi(t)= \sum_{j\in\N}\delta_{Y_j(t)}$ 
solves the ISDE 
\begin{equation}\label{ISDE:Ai}
\tag{Ai} 
dY_j(t)=dB_j(t)
+\lim_{r\to\infty}\left\{\sum_{\substack{k:\, k\not=j \\ |Y_k(t)|<r}}^{\infty}\frac{1}{Y_j(t)-Y_k(t)}-\int_{-r}^r \frac{\hat{\rho}(x)dx}{-x}\right\}dt,
\end{equation}
under $\mathbf{P}_{\Ai}^\xi$ with $\hat{\rho}(x)\equiv (\sqrt{-x}/\pi) {\bf 1}(x<0)$,
and  a labeled process $\mathbf{Z}^\nu =(Z_j^\nu)_{j\in\N}$ with 
$\Xi(t)= \sum_{j\in\N}\delta_{Z_j^\nu(t)}$ 
solves the ISDE 
\begin{equation} \label{ISDE:J}
\tag{$J_{\nu}$}
dZ_j^\nu(t)=2\sqrt{Z_j^\nu(t)}d{B}_j(t)+2(\nu+1)dt
+\sum_{\substack{k: \, k \not= j \\ k=1 }}^{\infty} 
\frac{4Z_j^\nu(t) \ dt}{Z_j^\nu(t)-Z_k^\nu(t)}.
\end{equation}
under $\mathbf{P}_{J_\nu}^\xi$ $\nu>-1$ \cite{o.isde, ot.1, ho.1}.

Set $\mathbf{P}_\star^{\mu_\star}= \int_{\mM} \mu_\star (d\xi) \mathbf{P}_\star^\xi$, 
for $\star \in \{\sin, \Ai\} \cup \{J_\nu\ ;\nu>-1\}$. $(\Xi(t), \mathbf{P}_\star^{\mu_\star})$ is then a reversible diffusion process with reversible measure $\mu_\star$. The second main result is the following theorem.

\begin{thm}\label{TH:SDE}
Let $\star \in \{\sin, \Ai\} \cup \{J_\nu\ ;\nu>-1\}$. 

\noindent {\rm (i)} A labeled process $\X =(X_j)_{j\in\N}$ associated with the determinantal process $(\Xi(t), \P_{\star})$ 
solves the ISDE ($\star$).

\noindent {\rm (ii)} Process $(\Xi(t), \P_{\star})$ coincides with process $(\Xi(t), \mathbf{P}_{\star}^{\mu_\star})$ in distribution. In particular, process $(\Xi(t), {\P}_{\star})$ is associated with the Dirichlet form $(\cE_{\star}, \cD_{\star})$.
\end{thm}

We remark that Tsai \cite{Tsai} recently constructed non-equilibrium dynamics for the Dyson model with $0<\beta<\infty$ by ISDE and obtained the result related to (i) in the above theorem.

Consider that the DPPs $\mu_{\sin}^N$, $\mu_{\Ai}^N$ and $\mu_{J_\nu}^N$ are probability measures on the configuration spaces of a finite number of particles defined by (\ref{mu_sin_N}), (\ref{mu_Ai_N}) and (\ref{mu_J_N}), respectively, and the bilinear forms $(\cE_\mu^a, \cD_\infty^\mu)$ with $(\mu, a(x))= (\mu_{\sin}^N, 1), (\mu_{\Ai}^N,1)$ and  $(\mu_{J_\nu}^N, 4x)$.
The closure $(\cE_{\star}^N, \cD_{\star}^N)$ of the bilinear form  $(\cE_{\mu_{\star}^N}^a, \cD_\infty^{\mu_\star^N})$ is a quasi-regular Dirichlet form, associated with the diffusion process
$(\Xi(t), \mathbf{P}_\star^{\xi^N})$, $\star \in \{\sin, \Ai\} \cup \{J_\nu\ ;\nu>-1\}$. We see that labeled processes associated with these diffusion processes solve SDEs (\ref{eqn:Dyson2}), (\ref{eqn:Airy2}) and (\ref{eqn:Bessel2}), respectively.
We take a labeled map $\mathfrak{l} : \mM \to \R^{\N} \oplus \bigoplus_{n=0}^\infty \R^n$ such that for $\xi\in\mM$
$$
|\mathfrak{l}_{j-1}(\xi)| \le |\mathfrak{l}_{j}(\xi)|, \quad j= 1,2,\dots, \xi(\R),
$$
and set 
\begin{equation}\label{:xil}
\xi_N^{\mathfrak{l}}=\sum_{j=1}^N \delta_{\mathfrak{l}_j(\xi)}
\quad \text{for } \xi\in \mM.
\end{equation}
Let $\X =(X_j)_{j\in\N}$ and $\X^N =(X_j^N)_{j=1}^{N}$ be the labeled processes associated with $(\Xi(t), \mathbf{P}_{\star}^{\xi})$ and $(\Xi(t), \mathbf{P}_{\star}^{\xi^N})$, respectively.
Note that $\X(0)=\mathfrak{l}(\xi)\equiv\x$ and $\X^N(0)=(\mathfrak{l}_j(\xi))_{j=1}^N\equiv\x^N$. 
We have then the following as a corollary of Theorem \ref{TH:SDE}.

\begin{cor}\label{C:23}
Let $\star \in \{\sin, \Ai\} \cup \{J_\nu\ ;\nu>-1\}$. 
{\rm (i)} 
 For $\mu_\star$ a.s. $\xi$,
$$
(\Xi(t),\mathbf{P}_{\star}^{\xi_N^{\mathfrak{l}}})\to (\Xi(t), \mathbf{P}_{\star}^{\xi}), \quad N\to\infty, 
$$
weakly on the path space $C([0,\infty),\mM)$.
 
\noindent {\rm (ii)} 
For $\mu_{\star}\circ \mathfrak{l}^{-1}$ a.s. $\x$, and $m\in\N$
$$
(X_1^N(t), X_2^N(t), \dots, X_m^N(t)) \to (X_1(t), X_2(t), \dots, X_m(t)), \quad N\to \infty,
$$
weakly on the path space $C([0,\infty),\R^m)$.
\end{cor}

Suppose that $\mu$ and $ \mu_N$ ($ N\in\mathbb{N}$) are probability measures on $\mM$, 
and $(\Xi(t), P)$ and $ (\Xi(t), P^N)$ are diffusion processes associated with the Dirichlet spaces given by  the closures of $(\cE_{\mu}^a, \cD_\infty^{\mu}, L^2(\mM ,\mu ))$ and 
$ (\cE_{\mu_N}^a, \cD_\infty^{\mu_{N}}, L^2(\mM ,\mu^N ))$, respectively. 
Let us consider the problem on the weak convergence of stationary processes. 
That is, 
$$
\mu_N \to \mu, \ N\to\infty \
\Rightarrow
(\Xi(t), P^N) \to (\Xi(t), P), \ N\to\infty
.$$
If the measures $ \mu_{N}$ and $ \mu $ are singular each other,  then such a convergence is not covered by a general theorem of convergence of diffusions associated with Dirichlet forms. 
We remark that Corollary \ref{C:23} (i) gives examples of such a convergence even if 
the measures $ \mu_{N}$ and $ \mu $ are singular each other. 
Recently, Kawamoto-Osada \cite{ko.1} also showed Corollary \ref{C:23} (ii) using a different method.

\SSC{Proof of theorems}
\subsection{Some properties of determinantal point processes}

Let $\mu$ be a probability measure on $\mM$
with correlation functions 
$\rho_m(\x_m)$, $\x_m \in \R^m$,
$m\in\N$.
Then, for $f\in\C_0(\R)$ and $\theta\in\R$
\begin{eqnarray}
\Psi(f,\theta)
\equiv \int_{\mM}\mu(d\xi) e^{\theta\int_{\R^d}f(x)\xi(dx)}
=\sum_{m=0}^\infty \frac{1}{m !}
\int_{\R^{m}} d\x_m
\prod_{j=1}^m \Big(e^{\theta f(x_j)}-1\Big)
\rho_m \Big( \x_m \Big),
\nonumber
\end{eqnarray}
where $d\x_m=\prod_{j=1}^{m} d x_{j}$.
Let ${\cal K}$ be a symmetric linear operator with kernel $K$.
In this section we assume that operator ${\cal K}$  satisfies 
{\bf Condition A} in \cite{ST03}.
Probability measure $\mu$ is called a DPP
with correlation kernel $K$,
if its correlation functions are given by
$$
\rho_m(\x_m)
=\det_{1 \leq j,k \leq m}
\Big[
K(x_j,x_k)
\Big].
$$
We often write $\rho$ for $\rho_1$, and $\rho(A)$ 
for $\int_A \rho(x) dx$, $A\in{\cal B}(\R)$.
The following lemmas are \cite[Lemma 3.1]{KT11b} and a modification of \cite[lemma 3.3]{KT11b}.

\begin{lem}
\label{L:31}
Let $\mu$ be a DPP. For any bounded closed interval $D$ of $\R$, we then have
\begin{equation}
\int_{\mM} \mu(d\eta) \Big| \eta(D)-\int_D \rho (x)dx\Big|^{2k}
\le ( 3 \rho(D) )^k, \quad k\in\N.
\nonumber 
\end{equation}
\end{lem}

\vskip 3mm

\begin{lem}
\label{L:32}
Let $\mu$ be a DPP and $(\Xi(t), P)$ be a stationary process with stationary measure $\mu$. If
\begin{equation}\label{l321}
\sum_{k\in\Z}|k|^{\ell}
\rho([g^{\kappa}(k),g^{\kappa}(k+1)])^m  <\infty,
\end{equation}
for some $\ell, m\in \N$, and $\kappa >0$,
then for $P$-a.s. $\Xi$ there exists $m_0=m_0(\Xi)\in\N$ such that
\begin{equation}\label{l322}
\Xi(t,[g^{\kappa}(k),g^{\kappa}(k+1)])\le m_0, \quad t=\frac{j}{|k|^\ell}, j=1,2,\dots, |k|^\ell, k\in\Z,
\end{equation}
where $g^{\kappa}$ is the function in (\ref{gk}). 
In particular
\begin{equation}\label{l323}
\lim_{m\to\infty}P\left(\Xi(t,[g^{\kappa}(k),g^{\kappa}(k+1)])\le m, \ t=\frac{j}{|k|^\ell}, j=1,2,\dots, |k|^\ell, k\in\Z\right)=1.
\end{equation}

\end{lem}

\noindent {\it Proof.} Since a DPP has the repulsive property, from condition (\ref{l321}) we derive
\begin{equation}
\nonumber 
\sum_{k\in\Z}|k|^{\ell}
\int_{[g^{\kappa}(k),g^{\kappa}(k+1)]^m}\rho_m(\x_m)d\x_m  <\infty,
\end{equation}
which implies
\begin{eqnarray}\nonumber
&&\sum_{k\in\Z}|k|^{\ell}\mu \Big(\xi(g^{\kappa}(k), g^{\kappa}(k+1))>m \Big) 
\\
&&=\sum_{k\in\Z}P\left(\Xi(t,[g^{\kappa}(k),g^{\kappa}(k+1)])> m, \mbox{ for some } t\in \left\{\frac{j}{|k|^\ell}; j=1,2,\dots, |k|^\ell\right\} \right)
<\infty.
\nonumber
\end{eqnarray}
By Borel-Cantelli's lemma we obtain (\ref{l321}), and then (\ref{l322}).
\qed

\subsection{Non-equilibrium dynamics}\label{S3.2}

In this subsection we review the results in \cite{KT09, KT10a, KT11a} on non-equilibrium dynamics of noncolliding diffusion processes.
For $\xi^N \in \mM$ with $\xi^N(\R)=N\in\N$ and $p 
\in \N_0 \equiv \N \cup \{0\}$,
we consider the product
$$
\Pi_{p}(\xi^N, w)= 
\prod_{x \in \supp \xi^{N}}
G \left( \frac{w}{x}, p \right)^{\xi^N(\{x\})},
\quad w \in \C,
$$
where 
\begin{equation}
G(u,p) = \left\{
   \begin{array}{ll}
\displaystyle{1-u},
& \mbox{if} \quad p=0  \\
& \\
\displaystyle{(1-u)\exp\left [ u+\frac{u^2}{2}+\cdots +\frac{u^p}{p} \right]},
& \mbox{if} \quad p\in\N.
   \end{array} \right. 
\nonumber
\end{equation}
The functions $G(u,p)$ are called
the Weierstrass primary factors \cite{Lev96}.
We then set
\begin{equation}
\Phi_p(\xi^N, z, w) 
\equiv \Pi_p(\tau_{-z}\xi^N \cap \{0\}^{\rm c}, w-z)
=\prod_{x\in\supp \xi^N \cap \{z\}^{\rm c}}
G\left(\frac{w-z}{x-z},p \right)^{\xi^N(\{x\})},
\quad w, z \in \C,
\nonumber
\end{equation}
where $\tau_z \xi(\cdot) \equiv \sum_{j \in \Lambda} \delta_{x_j+z}(\cdot)$ 
for $z \in \C$ and $\xi(\cdot) =\sum_{j \in \Lambda} \delta_{x_j}(\cdot)$.


Let $\u$ be the unlabeled map from $\cup_{n=0}^\infty \R^n$ to $\mM$ defined as $\u (x_1, \dots, x_n)= \sum_{j=1}^n \delta_{x_j}$. 
We consider the process $\u(\X(t))$, $t\in [0,\infty)$,
for the solution $\X(t)=(X_j(t))_{j\in\I_N}$ of (\ref{eqn:Dyson1}). We donote by $\P_{\sin}^{\xi^N}$ the distribution of the process starting from any fixed configuration $\xi^N$ with a finite number of particles.
It was proved in Proposition 2.1 of \cite{KT10a} that process $(\Xi(t), \P_{\sin}^{\xi^N})$ is determinantal with the correlation kernel $\mathbb{K}^{\xi^N}$
given by 
\begin{eqnarray}
\mathbb{K}_{\sin}^{\xi^N}(s, x; t, y)
&=& \frac{1}{2 \pi i} \int_{\R} du \,
\oint_{{\rm C}_{iu}(\xi^N)} dz \,  
p_{\sin}(s, x|z) 
\frac{\Pi_0(\tau_{-z}\xi^N, iu-z)}{iu-z}
p_{\sin}(-t, iu|y)
\nonumber\\
&&  - {\bf 1}(s > t)p_{\sin}(s-t, x|y),
\label{def:K_xi}
\end{eqnarray}
where ${\rm C}_w(\xi^N)$ denotes a closed contour on the
complex plane $\C$ encircling the points in 
$\supp \xi^N$ on the real line $\R$
once in the positive direction but not point $w$, 
and $p_{\sin}(t, x|y)$ is the generalized heat kernel:
\begin{equation}
p_{\sin}(t, x|y) 
= \frac{1}{\sqrt{2\pi |t|}}\exp\Big\{-\frac{(x-y)^2}{2t}\Big\}{\bf 1}(t\not= 0)
+ \delta(y-x){\bf 1}(t=0),
\quad t \in\R, \ x,y\in\C.
\nonumber
\end{equation}
In case $\xi^N$ is a configuration without any multiple points, i.e. $\xi^N \in \mM_0$,
(\ref{def:K_xi}) is rewritten as 
\begin{eqnarray}
\mathbb{K}_{\sin}^{\xi^N}(s, x; t, y)
&=& \int_{\R} \xi^N(dx') \, 
\int_{\R} d u \, p_{\sin}(s, x|x')
\Phi_0 (\xi^N, x',iu)p_{\sin}(-t, iu|y)
\nonumber\\
&& - {\bf 1}(s > t) p_{\sin}(s-t, x|y).
\nonumber
\end{eqnarray}

It was proved in \cite[Theorem 2.4]{KT10a} with  \cite[Theorem 1.4]{KT13} that
if $\xi\in\mX_{\sin}$, process 
$(\Xi(t), \P_{\sin}^{\xi_N^{\mathfrak{l}}})$
converges to the determinantal process with correlation kernel $\mathbb{K}^{\xi}$
as $N\to\infty$, weakly on $C([0,\infty),\mM)$.
In particular, when $\xi \in \mX_{\sin} \cap \mM_0$,
$\mathbb{K}_{\sin}^{\xi}$ is given by 
\begin{eqnarray}
\mathbb{K}_{\sin}^{\xi}(s, x; t, y)
&=& \int_{\R} \xi(dx') \, 
\int_{\R} du \, p_{\sin}(s, x|x')
\Phi_0 (\xi, x',iu)p_{\sin}(-t, iu|y)
\nonumber\\
&& - {\bf 1}(s > t) p_{\sin}(s-t, x|y).
\nonumber
\end{eqnarray}

\vskip 3mm


Let $\hat{\rho}^N, N\in\N$ be a sequence of nonnegative functions on $\R$ defined as
$$
\hat{\rho}^N(x)=\frac{\1_{(-4N^{2/3},0]}(x)}{\pi}\sqrt{-x\left(1+\frac{x}{4N^{2/3}}\right)}.
$$
Then, $\int_{\R} dx \ \hat{\rho}^N(x) =N$, $\int_{\R} dx \ \frac{\hat{\rho}^N(x)}{-x} =N^{1/3}$,
 and
$$
\hat{\rho}^N(x) \nearrow 
\hat{\rho}(x) = (\sqrt{-x}/\pi) {\bf 1}(x<0), \quad N \to \infty.
$$
Consider the process $\Y(t)=(Y_j(t))_{j\in\I_N}$ given by
\begin{equation}
Y_j(t) = X_j(t)+\frac{t^2}{4} +t \int_{\R}\frac{\hat{\rho}^N(x)dx}{x}, \quad
j\in\I_N,
\nonumber 
\end{equation}
with the solution
$\X(t)=(X_j(t))_{j\in\I_N}$ of (\ref{eqn:Dyson1}).
In other words, $\Y(t)$
satisfies the following SDE:
\begin{eqnarray}
&&dY_j(t) = dB_j(t)+\left( \frac{t}{2}-N^{1/3} \right) dt
+\sum_{\substack{k:\, k \not= j \\ k = 1 }}^{N}\frac{dt}{Y_j(t)-Y_k(t)},
\quad j\in\I_N.
\label{Ai_finite}
\end{eqnarray}
We denote by $\P_{\Ai}^{\xi^N}$ the distribution of the process $\u(\Y(t))$, $t\in [0,\infty)$
starting from any fixed configuration $\xi^N$ with a finite number of particles.
It was proved in Proposition 2.4 in \cite{KT09} that process $(\Xi(t),\P_{\Ai}^{\xi^N})$ is determinantal with correlation kernel $\mathbb{K}_{\hat{\rho}^N}^{\xi^N}$
given by 
\begin{eqnarray}
\mathbb{K}^{\xi^{N}}_{\hat{\rho}^N}(s,x ;t,y)
&=& \frac{1}{2\pi i}
\int_{\R} du \,
\oint_{{\rm C}_{iu}(\xi^N)} dz \,
q(0,s, x-z)
\frac{\Pi_0(\tau_{-z}\xi,iu-z)}{iu-z}
\nonumber\\
&&\qquad\times
\exp\bigg[(iu-z)\int_{\R}\frac{\hat{\rho}^N(v)dv}{v}\bigg]
q(t,0, iu-y)
\nonumber\\
&& - {\bf 1}(s>t)q(t, s, x-y),
\nonumber 
\end{eqnarray}
where $q(s,t,y-x), s, t \in \R, s \not= t, x, y \in \C$, is given by
$$
q(s, t, y-x) =
p_{\sin}\left(t-s, \left(y-\frac{t^2}{4} \right)-
\left( x-\frac{s^2}{4} \right) \right).
$$
Note that $q(s,t,y-x)$, $0\le s <t$, $x,y \in\R$ 
is the transition density function of process 
$B(t)+t^2/4$,
where $B(t), t \in [0, \infty)$ is the one-dimensional
standard BM.

Let $M_{\cA}(\xi)$ be the function defined as
\begin{equation}
M_{\cA}(\xi)
=\lim_{L\to\infty} \int_{0<|x|<L} \frac{\hat{\rho}(x)dx-\xi(dx)}{x}.
\nonumber
\end{equation}
For $\xi \in \mM$ with $M_{\cA}(\xi)<\infty$ and $z\in\C$, we define
\begin{eqnarray}
&&\Phi_{\Ai}(\xi,w) 
\equiv \exp \Bigg[ w M_{\cA}(\xi) \Bigg]
\Pi_1(\xi\cap \{0\}^{\rm c}, w), \quad w \in \C, 
\nonumber\\
&&\Phi_{\Ai} (\xi,z,w)\equiv \Phi_{\Ai}(\tau_{-z}\xi,w-z), \quad w, z \in \C.
\nonumber
\end{eqnarray}
We note that $\Phi_{\Ai} (\xi,z,w)$ exists finitely and $\Phi_{\cA}(\xi,z,w) \not \equiv0$, if $\xi\in\mX_{\Ai}$.

If $\xi\in \mX_{\Ai}$, sequence of the processes $(\Xi(t), \P_{\Ai}^{\xi_N^{\mathfrak{l}}})$
converges to the determinantal process $(\Xi(t), \P_{\Ai}^{\xi})$ with correlation kernel $\mathbb{K}^{\xi}_{\Ai}$ as $N\to\infty$, 
weakly on $C([0,\infty),\mM)$.
In particular, when 
$\xi \in \mX_{\Ai} \cap \mM_0$,
$\mathbb{K}^{\xi}_{\Ai}$ is given by 
\begin{eqnarray}
&& \mathbb{K}^{\xi}_{\Ai}(s,x;t,y)
=\int_{\R}\xi(dx') 
\int_{\R} du \, 
q(0,s, x-x') 
\Phi_{\Ai}(\xi, x', iu)
q(t, 0, iu-y)
\nonumber\\
&& \qquad\qquad\qquad - {\bf 1}(s>t)q(t, s, x-y).
\nonumber
\end{eqnarray}
This result was stated in \cite[Section 2.4]{KT09}, in which the weak convergence is in the sense of finite dimensional distributions. The convergence on $C([0,\infty),\mM)$ is derived using the same argument as \cite[Theorem 1.4]{KT13}. 


We consider a one-parameter family of $\mM^+$-valued processes 
with parameter $\nu > -1$, 
$\u({\bf Z}^{\nu}(t))$,
$t \in [0, \infty)$,
for the solution ${\bf Z}^{\nu}(t)=(Z_j^\nu (t))_{j\in\I_N}$ of SDE (\ref{eqn:Bessel1}). 
For a given configuration $\xi^N \in \mM^{+}$ with a finite number of particles, 
we denote by $\P_{J_\nu}^{\xi^N}$ the distribution of the process starting from $\xi^N$. 
In \cite[Theorem 2.1]{KT11a} it was proved that 
$(\Xi(t), \P_{J_\nu}^{\xi^N})$ 
is determinantal with correlation kernel 
\begin{eqnarray}
\mathbb{K}^{\xi^N}_{J_\nu}(s, x; t, y)
&=& \frac{1}{2 \pi i} 
\int_{-\infty}^{0} du
\oint_{{\rm C}_u(\xi^N)} dz  \,  
p^{(\nu)}(s, x|z) \frac{\Pi_0(\tau_{-z}\xi^N, u-z)}{u-z}
p^{(\nu)}(-t, u|y)
\nonumber\\
&& - {\bf 1}(s > t)p^{(\nu)}(s-t, x|y),
\label{eqn:KN1a}
\end{eqnarray}
where for $t \in \R$ and $x, y \in \C$
\begin{eqnarray}
p^{(\nu)}(t, y|x)
&=& \frac{1}{2|t|} \left(\frac{y}{x} \right)^{\nu/2}
\exp \left(-\frac{x+y}{2t} \right)
I_{\nu}\left(\frac{\sqrt{xy}}{|t|} \right)
{\bf 1}(t \not=0, x\not=0)
\nonumber\\
&+&\frac{y^\nu}{(2|t|)^{\nu+1}\Gamma(\nu+1)}
\exp \left(-\frac{y}{2t} \right){\bf 1}(t \not=0, x=0)
+\delta(y-x) {\bf 1}(t=0).
\label{eqn:pnu-}
\end{eqnarray}
Note that for $t\ge 0$ and $x,y\in\R_+$, $p^{(\nu)}(t, y|x)$ is the transition density function of a $2(\nu+1)$-dimensional squared Bessel process.
If $\xi\in\mX_{J_\nu}$,  process $(\Xi(t), \P_{J_\nu}^{\xi_N^{\mathfrak{l}}} )$ converges to the determinantal process 
with correlation kernel $\mathbb{K}^{\xi}_{\nu}$
as $N\to\infty$ weakly on $C([0,\infty),\mM))$.
In particular, when 
$\xi \in \mX_{J_\nu} \cap \mM_0$,
$\mathbb{K}^{\xi}_{J_\nu}$ is given by 
\begin{eqnarray}
\mathbb{K}^{\xi}_{J_{\nu}}(s, x; t, y)
&=&  
\int_{0}^{\infty} \xi(dx') 
\int_{-\infty}^{0} du \,
 p^{(\nu)}(s, x|x') \Phi_0(\xi, x', u)  p^{(\nu)}(-t, u|y)
\nonumber\\
&& - {\bf 1}(s > t)p^{(\nu)}(s-t, x|y).
\nonumber
\end{eqnarray}
The result was proved in \cite[Theorem 2.2]{KT11a} in which the weak convergence is in the sense of finite dimensional distributions.
The weak convergence on $C([0,\infty),\mM)$ can be proved from the tightness of $\{ \P_{J_\nu}^{\xi_N^{\mathfrak{l}}}\}_{N\in\N}$, which is derived by the same procedure to show \cite[Theorem 1.4]{KT13} using determinatal martingales representation in \cite[Theorem 1.2]{K14} instead of complex Brownian motion representation in \cite[Theorem 1.1]{KT13}.

The following proposition is a combination of the results in \cite{KT09, KT10a, KT11a, KT11b}  and the tightness of $\{\P_\star^{\xi_n}\}_{n\in\N}$, $\star \in \{\sin, \Ai\} \cup \{J_\nu\ ;\nu>-1\}$, which is derived by the same argument as that of $\{ \P_{\star}^{\xi_N^{\mathfrak{l}}}\}_{N\in\N}$.

\begin{prop} \label{Prop:33}
 Let $\star \in \{\sin, \Ai\} \cup \{J_\nu\ ;\nu>-1\}$.

\noindent {\rm (i)}  
Suppose that $\xi, \xi_n \in \mX_{\star}, n\in\N$.
If $\xi_n$ converges to $\xi$ in $\mX_{\star}$, then
\begin{equation}\label{conv:xi}
(\Xi(t),\P_\star^{\xi_n})
\to (\Xi(t),\P_\star^{\xi}), \quad n \to \infty
\end{equation}
weakly on $C([0,\infty),\mM)$. 

\noindent {\rm (ii)}
Determinantal process $(\Xi(t), \P_\star)$ with the correlation kernel $\mathbb{K}_\star$ is identical in distribution to process $(\Xi(t), \P_\star^{\mu_\star})$, where $\P_\star^\nu =\int_{\mM}\P_\star^\xi \ \nu(d\xi)$ for probability measure $\nu$ on $\mM$.
\end{prop}

\vskip 3mm

\subsection{Path property of noncolliding processes}

Let $\X(t)=(X_1(t), X_2(t), \dots, X_N(t))$ be the Dyson model defined by (\ref{eqn:Dyson1}). We introduce the following version of the Dyson model 
\begin{equation}\nonumber
\tilde{X}_j(t)=e^{-\gamma_N t}X_j(\tau_N(t)),
\end{equation}
where $\gamma_N=\frac{1}{2N}$ and $\tau_N (t)=(e^{2t\gamma_N }-1)/2\gamma_N$. 
Process $\tilde{\X}(t)=(\tilde{X}_j(t))_{j=1,\dots,N}$ then solves the SDE
\begin{equation}
d\tilde{X}_j(t)=d\tilde{B}_j(t)-\frac{\tilde{X}_j(t)dt}{2N} + \sum_{\substack{k:\, k \not= j \\ k = 1 }}^{N}\frac{dt}{\tilde{X}_j(t)-\tilde{X}_k(t)}, \quad j\in\I_N,
\label{eqn:Dyson2}
\end{equation}
where the $\tilde{B}_j(t)$'s are independent one-dimensional standard BMs. The eigenvalue distribution in GUE,
\begin{equation}\label{GUE_N}
m_{\rm GUE}^N(d\x_N)=\frac{1}{Z}\prod_{1\le i<j\le N}|x_i-x_j|^2
\exp\left\{-\frac{1}{2N}\sum_{k=1}^N x_k^2\right\}d\x_N
\end{equation}
is a reversible probability measure of process $\tilde{\X}(t)$. Hereafter, the notation $Z$ is used to denote the normalizing constant. 
We denote the distribution of $\u(\X(t))$ with $\u(\X(0))=\xi^N$ by $\tilde{\P}_{\sin}^{\xi^N}$. Process $(\Xi(t), \tilde{\P}_{\sin}^{\xi^N})$ has the probability measure 
\begin{equation}\label{mu_sin_N}
\mu_{\sin}^N \equiv m_{\rm GUE}^N\circ \u^{-1}
\end{equation}
as the reversible measure. 
Measure $\mu_{\sin}^N$ is the DPP with the correlation kernel
\begin{equation}
\mathbb{K}_N(x,y) = 
\frac{1}{\sqrt{2N}}
\sum_{k=0}^{N-1}
\varphi_{k}\left(\frac{ x}{\sqrt{2N}}\right)
\varphi_{k}\left(\frac{ y}{\sqrt{2N}}\right),
\nonumber 
\end{equation}
where $h_k=\sqrt{\pi}2^k k!$ and 
\begin{equation}
\varphi_k(x)=\frac{1}{\sqrt{h_k}}e^{-x^2/2}H_k(x),
\nonumber 
\end{equation}
with Hermite polynomials $H_k, k\in \N_0$.

We also introduce process $\tilde{\Y}(t)=(\tilde{Y}_j(t))_{j\in\I_N}$ obtained from the Dyson model in (\ref{eqn:Dyson2}) by the transformation given by $N^{-1/3}\tilde{X_j}(N^{2/3}t)-2N^{2/3}$, which solves the SDE 
\begin{equation}
d\tilde{Y}_j(t)=dB_j(t)-\frac{\tilde{Y}_j(t)+2N^{2/3}}{2N^{1/3}}dt + 
\sum_{\substack{k:\, k \not= j \\ k = 1 }}^{N} 
\frac{dt}{\tilde{Y}_j(t)-\tilde{Y}_k(t)}, \quad j\in\I_N.
\label{eqn:Airy2}
\end{equation}
We denote the distribution of process $\u(\tilde{\Y}(t))$ with $\u(\tilde{\Y}(0))=\xi^N$ by $\tilde{\P}_{\Ai}^{\xi^N}$. 
Process $(\Xi(t), \tilde{\P}_{\Ai}^{\xi^N})$
has
\begin{equation}\label{mu_Ai_N}
\mu_{\Ai}^N\equiv m_{\rm Ai}^N\circ \u^{-1}
\end{equation}
with 
\begin{equation}\nonumber
m_{\rm Ai}^N(d\x_N)=\frac{1}{Z}\prod_{1\le i<j\le N}|x_i-x_j|^2
\exp\left\{-\frac{1}{2}\sum_{k=1}^N (x_k-N^{1/3})^2\right\}d\x_N
\end{equation}
as the reversible measure.
For the noncolliding squared Bessel process in (\ref{eqn:Bessel1}) we set $\tilde{X}_j^\nu(t)=e^{-\gamma_N t}X_j^\nu(\tau_N(t))$. Process $\tilde{\X}^\nu(t)=(\tilde{X}_j^\nu(t))_{j\in\I_N}$ then solves the SDE
\begin{equation}
d\tilde{X}_j^\nu(t)=2\sqrt{\tilde{X}_j^\nu(t)}d{B}_j(t)-\frac{\tilde{X}_j^\nu(t)dt}{N}+2(\nu+1)dt
+ \sum_{\substack{k:\, k \not= j \\ k = 1 }}^{N} \frac{4\tilde{X}_j^\nu(t) \ dt}{\tilde{X}_j^\nu(t)-\tilde{X}_k^\nu(t)}.
\label{eqn:Bessel2}
\end{equation}
The eigenvalue distribution in chGUE,
\begin{equation}\label{chGUE_N}
m_{J_\nu}^N(d\x_N)=\frac{1}{Z}\prod_{1\le i<j\le N}|x_i-x_j|^2 \prod_{j=1}^N x_j^{\nu+1/2}
\exp\left\{-\frac{1}{2}\sum_{k=1}^N x_k\right\}d\x_N,
\end{equation}
is a reversible probability measure of $\tilde{\X}^\nu(t)$.
We denote the distribution of process 
$\u(\tilde{\X}^\nu(t))$ by $\tilde{\P}_{J_\nu}^{\xi^N}$. Process $(\Xi(t), \tilde{\P}_{J_\nu}^{\xi^N})$ has the probability measure 
\begin{equation}\label{mu_J_N}
\mu_{J_\nu}^N \equiv m_{J_\nu}^N\circ \u^{-1}
\end{equation} 
as the reversible measure.
Measure $\mu_{J_\nu}^N$ is
the DPP with correlation kernel 
$$
\mathbb{K}^{(\nu)}_{N}(x,y)= 
\frac{1}{2N}\sum_{k=0}^{N-1}
\varphi_k^\nu \left(\frac{x}{2N}\right)
\varphi_k^\nu \left(\frac{y}{2N}\right),
$$
where 
$$
\varphi_k^\nu (x) = 
\sqrt{\frac{\Gamma(k+1)}{\Gamma(\nu+k+1)}}x^{\nu/2}
L_k^{\nu}(x)e^{-x/2},
$$ 
with the Laguerre polynomials $L_k^{\nu}(x), k\in\N_0$
with parameter $\nu>-1$ \cite{KT07a}.
From the construction of $\tilde{\P}_{\star}^{\xi^N}$, $\star \in \{\sin, \Ai\} \cup \{J_\nu\ ;\nu>-1\}$ the following is readily derived from Proposition \ref{Prop:33}.

\begin{prop}
\label{Prop:34}
Let $\star \in \{\sin, \Ai\} \cup \{J_\nu\ ;\nu>-1\}$.
Suppose that $\xi\in \mX_{\star}$. Then
$$
(\Xi(t), \tilde{\P}_{\star}^{\xi_N^{\mathfrak{l}}}) \to (\Xi(t), {\P}_{\star}^{\xi}), \quad N\to\infty,
$$
weakly on $C([0,\infty),\mM)$, where $\xi_N^{\mathfrak{l}}$ is the configuration with a finite number of particles given in (\ref{:xil}).
\end{prop}

\noindent {\bf Remark}
\; (i) Let $\star \in \{\sin, \Ai\} \cup \{J_\nu\ ;\nu>-1\}$. The reversible diffusion process $(\Xi(t), \tilde{\P}_{\star}^{\xi^N})$ is associated with the Dirichlet form $(\cE_\star^N, \cD_\star^N)$ introduced in Section 2.

\noindent (ii) We set $V_j^{\nu}(t)= \sqrt{\tilde{Z}_j^\nu(t)}$.
Process ${\bf V}(t)=(V_j(t))_{j\in\I_n}$ then solves the SDE
\begin{eqnarray}
&&dV_j^\nu(t)=d{B}_j(t)-\frac{V_j(t)dt}{2N}+\frac{2\nu+1}{2V_j(t)}dt
+ \sum_{\substack{k:\, k \not= j \\ k = 1 }}^{N} \frac{2V_j (t)\ dt}{V_j(t)^2-V_k(t)^2},
\label{sde:Bessel}
\end{eqnarray}
where we impose a reflecting wall at the origin for $\nu \in (-1,0)$.

\noindent (iii) Both processes $(\Xi(t), {\P}_{\Ai}^{\xi_N^{\mathfrak{l}}})$ and $(\Xi(t), \tilde{\P}_{\Ai}^{\xi_N^{\mathfrak{l}}})$ converge to process $(\Xi(t), {\P}_{\Ai}^{\xi})$ as $N\to\infty$ with the vague topology. However, when considering the limits of the associated SDEs (\ref{Ai_finite}) and (\ref{eqn:Airy2}), heuristically, they differ by $\frac{t}{2}dt$.
The difference stems from the fact that process $(\Xi(t), {\P}_{\Ai}^{\xi^N})$ is not reversible and particles far from the origin move in strong drifts. Then process $(\Xi(t), {\P}_{\Ai}^{\xi_N^{\mathfrak{l}}})$ does not converge to $(\Xi(t), {\P}_{\Ai}^{\xi})$ as $N\to\infty$ with a stronger topology that preserves the Markov property.
This would be related to the relation between the uniqueness of the strong solution of ISDE and the triviality of the tail $\sigma$-field of the path space \cite{ot.2}.
\vskip 3mm

Under $\P_\star$, $\star = \{\sin, \Ai\}\cup \{J_\nu, \nu>-1\}$, $\Xi\in C([0,\infty),\mM)$ is a non-explosion and non-crossing path and can be naturally lifted to the labeled path $\X=(X_j)_{j\in\N}$ such that $\Xi(t)=\sum_{j\in\N}\delta_{X_j(t)}$ and $X_j\in C([0,\infty), \R$), $j\in\N$. We write $X\in \Xi$ if $X\in \{X_j\}_{j\in\N}$.

\begin{lem}\label{L:35}
{\rm (i)} Let $\star\in \{\sin, \Ai\}$.
For each $T>0$, there exists a positive constant $C_\star$ such that
for any interval $D$ of $\R$ and $\varepsilon >0$
\begin{equation}
\P_{\star}\left(
\exists X \in \Xi \mbox{ s.t. }
X(0)\in D, \ 
\sup_{t\in [0,T]}\left|X(t)-X(0) \right|>\varepsilon\right)
\le C_\star(\rho_{\star}(D)\vee 1)\, {\rm Erf}\left(\frac{\varepsilon}{\sqrt{T}}\right) , 
\nonumber
\end{equation}
where  ${\rm Erf}(a)=\int_a^\infty \frac{1}{\sqrt{2\pi}}e^{-x^2/2}dx$.

\noindent{\rm (ii)}
Let $\star \in \{J_{\nu}, \nu>-1\}$.
For each $T>0$, there exists a positive constant $C_\star$ such that
for any interval $D$ of $[0,\infty)$ with $|D|\ge 1$ and $\varepsilon >0$
\begin{equation}
\P_{\star}\left(
\exists X \in \Xi \mbox{ s.t. }
X(0)\in D, \  
\sup_{t\in [0,T]}\left|\sqrt{X(t)}-\sqrt{X(0)} \right|>\varepsilon\right)
\le C_\star(\rho_{\star}(D)\vee 1) \, {\rm Erf}\left(\frac{\varepsilon}{\sqrt{T}}\right) .
\nonumber
\end{equation}
\end{lem}

To prove Lemma \ref{L:35} we prepare a lemma concerning on reversible diffusions 
in $ \mathbb{R}^N $. 

\begin{lem} \label{L:36}
Let $\Y(t)=(Y_j(t))_{j\in\I_N}$ be a reversible diffusion process satisfying 
$$
Y_j(t)=Y_j(0)+B_j(t)+\int_0^tb_j(\Y(s))ds,
\quad j\in\I_N, \ t\in [0,T],
$$
with some measurable function $b=(b_1,\dots,b_N)$. 
Here the $B_j(t), 1\le j \le N$ are independent BMs, and $\Y(0)=(Y_j(0))_{j\in\I_N}$ is a random variable distributed by the reversible probability measure, which is independent of the BMs. 
Putting $\rho_{Y}(A)=\sum_{j=1}^NP(Y_j(0)\in A)$ for $A\in {\cal B}(\R)$, and
\begin{align}\label{:36a}&
c=\sum_{\ell=0}^\infty 
\frac{\rho_{Y}(A_\ell)}{\rho_Y(D) \vee 1}
\frac{{\rm Erf}(\{2(\ell\vee 1)-1\}\varepsilon/\sqrt{T})}{{\rm Erf}(\varepsilon/\sqrt{T})},
\end{align}
we then obtain 
\begin{align} & 
P \left(Y_j(0)\in D, \, 
\sup_{t\in [0,T]}\left|Y_j(t)-Y_j(0) \right|>\varepsilon \mbox{ for some $j$}\right)
 \le 4(1+c)(\rho_Y(D)\vee 1) \, {\rm Erf}\left(\frac{\varepsilon}{\sqrt{T}}\right) . \label{est:LZ}
\end{align}\end{lem}

\noindent {\it Proof.} 
We use a Lyons-Zheng decomposition for the proof (see for example Section 5.7 in \cite{FOT}).
We deduce from the Lyons-Zheng decomposition that $\hat{\Y}(t)=\Y (T-t)$ satisfies
$$
\hat{Y}_j(t)=\hat{Y}_j(0)+\hat{B}_j(t)+\int_0^tb_j(\hat{\Y}(s))ds,
\quad j\in\I_N, \ t\in [0,T],
$$
with independent BMs $\hat{B}_j(t)$, $1\le j\le N$. Then
\begin{eqnarray}\nonumber
Y_j(t)-Y_j(0)&=& \hat{Y}(T-t)-\hat{Y}(T)=\hat{B}_j(T-t)-\hat{B}_j(T)-\int_0^tb_j(\Y(s))ds
\\
&=&\frac{1}{2}\left(B_j(t)+ \hat{B}_j(T-t)-\hat{B}_j(T)\right).
\label{LZ}
\end{eqnarray}
Putting $A_\ell=A_\ell(\varepsilon)=\{x\in\R : \inf_{y\in D}|x-y| \in [\ell\varepsilon,(\ell+1)\varepsilon) \}$, $\ell\in \N \cup \{0\}$, we have 
\begin{eqnarray}\nonumber
&&P\left(Y_j(0)\in D, 
\sup_{t\in [0,T]}\left|B_j(t)\right|\le \varepsilon,
\sup_{t\in [0,T]}\left|Y_j(t)-Y_j(0) \right|>\varepsilon \right)
\nonumber\\
&&\le \sum_{\ell=0}^\infty  
P\left(
Y_j(0)\in D, \, 
\sup_{t\in [0,T]}\left|B_j(t)\right|\le \varepsilon, \, 
Y_j(T)\in A_\ell, \, 
\sup_{t\in [0,T]}\left|Y_j(t)-Y_j(0) \right|\ge (\ell \vee 1)\varepsilon\right) 
\nonumber
\\
&&\le \sum_{\ell=0}^\infty  
P\left(
Y_j(0)\in D, \, Y_j(T)\in A_\ell, \, 
\sup_{t\in [0,T]}\left|\hat{B}_j(T-t)-\hat{B}_j(T) \right|\ge \{2(\ell\vee 1)-1\}\varepsilon\right) 
\ \text{ by (\ref{LZ})}
\nonumber 
\\
&&\le \sum_{\ell=0}^\infty  
P\left(
Y_j(T)\in A_\ell, \
\sup_{t\in [0,T]}
\left|\hat{B}_j(T-t)-\hat{B}_j(T) \right|\ge \{2(\ell\vee 1)-1\}\varepsilon \right).
\nonumber
\end{eqnarray}
Then from this we deduce that 
\begin{eqnarray}\nonumber
&&P\left(Y_j(0)\in D, \ 
\sup_{t\in [0,T]}\left|Y_j(t)-Y_j(0) \right|>\varepsilon
\ \mbox{ for some $j$}\right)
\\ \nonumber
&&\le \sum_{j=1}^N
P\left(Y_j(0)\in D, \ 
\sup_{t\in [0,T]}\left|Y_j(t)-Y_j(0) \right|>\varepsilon\right)
\nonumber\\
&&\le \sum_{j=1}^N P\left(
Y_j(0)\in D, \ 
\sup_{t\in [0,T]}\left|B_j(t)\right|>\varepsilon\right)
\nonumber\\
&&\quad \quad +\sum_{\ell=0}^\infty \sum_{j=1}^N P\left(Y_j(T)\in A_\ell, \
\sup_{t\in [0,T]}\left|\hat{B}_j(T-t)-\hat{B}_j(T) \right|\ge \{2(\ell\vee 1)-1\}\varepsilon\right).
\nonumber
\end{eqnarray}
Note that $\hat{\Y}(0)$ is independent of $\hat{B}_j(t)$, $ 1\le j\le N$, since $\Y(0)$ is independent of $B(t)$, $1\le j\le N$. We then have that the right-hand side of the above inequality is bounded by 
\begin{eqnarray}
&&\sum_{j=1}^N P\left(Y_j(0)\in D\right)
P\left(\sup_{t\in [0,T]}\left|B_j(t)\right|>\varepsilon\right)
\nonumber\\
&&\quad +\sum_{\ell=0}^\infty \sum_{j=1}^N P\left(Y_j(T)\in A_\ell\right)
P\left(\sup_{t\in [0,T]}\left|\hat{B}_j(T-t)-\hat{B}_j(T) \right|\ge \{2(\ell\vee 1)-1\}\varepsilon \right)
\nonumber\\
&&\le 4 \sum_{j=1}^N P\left(Y_j(0)\in D\right) \, {\rm Erf}\left(\frac{\varepsilon}{\sqrt{T}}\right) 
+4\sum_{\ell=0}^\infty \sum_{j=1}^NP(Y_j(0)\in A_\ell )
{\rm Erf}\left(\frac{\{2(\ell\vee 1)-1\}\varepsilon}{\sqrt{T}}\right),
\nonumber
\\
&&= 4 \rho_Y(D) \, {\rm Erf}\left(\frac{\varepsilon}{\sqrt{T}}\right) 
+4\sum_{\ell=0}^\infty \rho_{Y}(A_\ell) 
{\rm Erf}\left(\frac{\{2(\ell\vee 1)-1\}\varepsilon}{\sqrt{T}}\right). 
\nonumber
\end{eqnarray}
Here we used the estimate $P\left(\sup_{t\in [0,T]}\left|B_j(t)\right|>\varepsilon\right)
\le 4{\rm Erf}(\varepsilon/\sqrt{T})$ and the reversibility of the process. 
 (\ref{est:LZ}) follows from this and (\ref{:36a}) immediately. 
\qed 

\bigskip 

\noindent 
{\em Proof of Lemma \ref{L:35}.} 
We apply Lemma \ref{L:35} to SDEs (\ref{eqn:Dyson2}), (\ref{eqn:Airy2}), and (\ref{sde:Bessel}).
By simple calculation, for the solutions of these SDEs, $c$ is a constant independent of $D$ and $\varepsilon$.
We then see that there exists a positive constant $C_\star' >0$ such that for $\star\in\{\sin, \Ai\}$
\begin{equation}
\tilde{\P}_{\star}^{\mu_{\star}^N}\left(
\exists X \in \Xi_{\star} \mbox{ s.t. }
X(0)\in D, \ 
\sup_{t\in [0,T]}\left|X(t)-X(0) \right|>\varepsilon\right)
\le C_\star' (\rho_{\star}^N(D)\vee 1) \, {\rm Erf}\left(\frac{\varepsilon}{\sqrt{T}}\right) , \nonumber
\end{equation}
and for $\star \in \{
J_{\nu}, \nu>-1\}$
\begin{equation}
\tilde{\P}_{\star}^{\mu_{\star}^N}\left(
\exists X \in \Xi_{\star} \mbox{ s.t. }
X(0)\in D, \ 
\sup_{t\in [0,T]}\left|\sqrt{X(t)}-\sqrt{X(0)} \right|>\varepsilon\right)
\le C_\star' (\rho_{\star}^N(D)\vee 1) \, {\rm Erf}\left(\frac{\varepsilon}{\sqrt{T}}\right) , \nonumber
\end{equation}
where $\rho_{\star}^N(D)=\int_D \rho_{\star}^N(x)dx$ with  the density (the first correlation function) $\rho_{\star}^N(x)$ of $\mu_{\star}^N$.
Since $\tilde{\P}_{\star}^{\mu_{\star}^N}$ converges to $\P_{\star}$ weakly on $C([0,\infty),\mM)$, as $N\to\infty$ (see for instance \cite[Section 7]{KT07a}), we obtain the desire result by simple observation. \qed

\subsection{Proof of Theorem \ref{TH:SM}}

Let $\mX$ be a subset of $\mM$ and $(\Xi(t), P_{\xi})$, $\xi\in \mX$ be a stochastic process. Suppose that
\begin{align}\label{:34a}&
P_{\xi_n}(\cdot) \to P_{\xi}(\cdot), \quad \mbox{ if $\xi_n \to \xi $ in $\mX$}.
\end{align}
Further suppose
$\zeta\in\mX$ satisfying
$P_{\zeta}(\Xi(t) \mbox{ is continuous in $\mX$})=1$. 
If $(\Xi(t), P_{\zeta})$ is Markovian, then it is strong Markovian. 

Let $\star \in \{\sin, \Ai\} \cup \{J_\nu\ ;\nu>-1\}$. 
From Proposition \ref{Prop:33} \thetag{1} 
we see that $(\Xi(t), \P_{\star})$ satisfies \eqref{:34a}.  
Hence from the Markov property of process $(\Xi(t), \P_{\star})$  given in \cite{KT11b} and 
Proposition \ref{Prop:33} (ii), we deduce 
Theorem \ref{TH:SM} from the following proposition.

\vskip 3mm

\begin{prop}\label{Prop:3.7}
Let $\star \in \{\sin, \Ai\} \cup \{J_\nu\ ;\nu>-1\}$. Then process $(\Xi(t), \P_{\star})$ is continuous in $\mX_{\star}$.
\end{prop}

\noindent Proof. This proposition is derived from the following claims:

\noindent 1) $\Xi(t)$ has a vaguely continuous path.

\noindent 2) For each $T\in\N$, there exist $\varepsilon\in (0,1), \kappa\in (1/2, \kappa_\star)$ such that
\begin{equation}\nonumber
\lim_{L\to\infty}\lim_{m\to\infty}\P_{\star}\left(
\Xi(t) \in \mX_{L,m}^{\varepsilon, \kappa},
\ t\in [0,T] \right)=1.
\end{equation}

\noindent Claim 1) has already shown in Section \ref{S3.2}. We remark that it is also derived by Kolmogorov's criterion:
for any polynomial function $f$, which is a smooth function on $\mM$ givne in (\ref{pol}),
$$
\E_{\star}[|f(\Xi(t))-f(\Xi(s))|^\beta]\le C|t-s|^\alpha,
\quad 0\le s<t \le T< \infty,
$$
for some $\alpha>1, \beta> 0$, and $C>0$, which is readily proved from Lemma \ref{L:35}.

Claim 2) is derived from two estimates
\begin{equation}\label{claim2_1}
\lim_{L\to\infty}\P_{\star}\left(\left|\rho_{\star}(D_L) - \Xi(t,D_L)\right|\le L^\varepsilon, \ t\in [0,T] \right)=1, 
\ D_L=[0,L], [-L,0],
\end{equation}
and
\begin{equation}\label{claim2_2}
\lim_{m\to\infty}\P_{\star}\left(
\Xi(t,[g^\kappa(k),g^\kappa(k+1)])\le m,
\ t\in [0,T] \ k\in\Z \right)=1.
\end{equation}
From Lemma \ref{L:31}, there exist $m'\in\N$ and $p<m'-1$
such that
\begin{eqnarray}
\label{eqn:74}
&&\int_{\mM} \mu_{\star}(d\xi) \Big| \rho_{\star} (D_L)- \xi(D_L)\Big|^{m'}
= {\cal O}(L^p), \quad L\to\infty.
\end{eqnarray}
Taking $\varepsilon \in ((p+1)/m',1)$ and using Chebyshev's inequality with (\ref{eqn:74}),
we can find a positive constant $C$ such that
$$
\mu_{\star} \left(|\rho_{\star} (D_L)-\xi(D_L)|\ge L^{\varepsilon} \right)
\le CL^{p-m'\varepsilon}.
$$
Since $p-m'\varepsilon <-1$, we have
$$
\sum_{L=1}^\infty\mu \left(|\rho_{\star} (D_L)-\xi(D_L)| )
|\ge L^{\varepsilon} \right)<\infty.
$$
By Borel-Cantelli's lemma and the stationarity of the process, for any $k\in \N$
$$
\lim_{L\to\infty}\P_{\star}\left(\left|\rho_{\star}(D_L) - \Xi(t,D_L)\right|\le L^\varepsilon, \ t=\frac{j}{k}, \ j=1,2,\dots, kT \right)=1.
$$
(\ref{claim2_1}) is then derived from Lemma \ref{L:35}. Estimate (\ref{claim2_2}) is derived from Lemmas \ref{L:32} and \ref{L:35}. 
In fact Lemma \ref{L:35} (i) implies that for $\star\in\{\sin, \Ai \}$ and $k\in\Z$,
\begin{eqnarray}
&&\P_{\star}\left(
\exists X \in \Xi \mbox{ s.t. }
X(0)\in [g^{\kappa}(k),g^\kappa(k+1)], 
\sup_{u\in [0,1/k^{\ell}]}\left|X(u)-X(0) \right|> |k|^{\kappa-1}\right)
\nonumber\\
&&
\le C_\star\{\rho_{\star}([g^{\kappa}(k),g^\kappa(k+1)])\vee 1\}
{\rm Erf}(|k|^{\kappa+\ell/2-1}).
\nonumber
\end{eqnarray}
Take $\ell\ge 2$ and choose $m\in\N$ such that (\ref{l321}) holds.
Using Borel-Cantelli's lemma with simple calculations,
we see that 
for $\P_\star$-a.s. $\Xi$, there exists $k_0=k_0(\Xi)\in\N$ such that for any $k\in\Z$ with $|k|\ge k_0$,
if $X \in \Xi$ and $s\in\{\frac{j}{|k|^\ell}; j=1,2,\dots, |k|^\ell T \}$ satisfy
$$
X(s)\in [g^{\kappa}(k),g^\kappa(k+1)],
$$ 
then
$$
\sup_{u\in [0,1/k^{\ell}]}\left|X(s+u)-X(s) \right|\le |k|^{\kappa-1}.
$$
Combining this fact with (\ref{l322}),
we see that for $\P_\star$-a.s. $\Xi$ there exists $m_1=m_1(\Xi)\in\N$ such that
\begin{equation}
\Xi(t,[g^{\kappa}(k),g^{\kappa}(k+1)])\le m_1, \ , t\in [0,T], \ k\in\Z.
\end{equation}
We thus obtain (\ref{claim2_2}).
For $\star\in\{j_\nu, \nu>1\}$ we can obtain the desire result by the same argument as above with Lemma \ref{L:35} (ii).
This completes the proof.
\qed

\subsection{Proof of Theorem \ref{TH:SDE}}

We now introduce Dirichlet forms describing 
$ k $-labeled dynamics. 
For this we recall the definition of Palm and Campbell measures. 
Let $  \x_k =(x_1,\ldots,x_k)\in \R^{k}$. 
We set 
$$
\mu _{ \x_k  } = 
\mu (\cdot - \sum_{i=1}^{k} \delta _{x_i} | \ 
\xi( x_i )\ge 1 \text{ for }i=1,\ldots,k),
$$
which is called the (reduced) Palm measure of $\mu$.
The Campbell measure for probability measure $\mu$ is then given by
$$
\nu^{k}(d\x_k d\eta)= \mu_{\x_k}(d\eta)\rho_k(\x_k)d\x_k.
$$
Here $\rho _k : \R^k \to [0,\infty)$ 
is the $ k $-th correlation function of $ \mu $ 
 and $d \x_k =dx_1\cdots dx_k$ 
is the Lebesgue measure on $ \R^{k} $ as before. 
Let $\mathcal{D}^k_0 = C_0^\infty (\R^k)\otimes\mathcal{D}_\infty$.
For $f,g \in \mathcal{D}_0^k$, let $\nabla^k[f,g]$ be such that
\begin{equation}
\label{nabla}
\nabla^k[f,g](\x_k,\xi)=\frac{1}{2}\sum_{j=1}^k a(x_j)\frac{\partial}{\partial x_j}f(\x_k,\xi)\frac{\partial}{\partial x_j}g(\x_k,\xi).
\end{equation}
We set $\D^{a,k}$ as
\begin{equation}\nonumber
\D^{a,k}[f,g][\x_k,\xi]=\nabla^k[f,g](\x_k,\xi)+\D [f(\x_k,\cdot),g(\x_k,\cdot)](\xi).
\end{equation}
We consider the bilinear form $(\cE_{\nu^k}^{a,k} , {\cal D}_\infty^{\nu^k})$ on $L^2(\R^k\times\mM,\nu^k)$ defined as
\begin{eqnarray}
&&\cE_{\nu^k}^{a,k} (f,g)= \int_{\mM}{\Bbb D}^{a,k}(f,g)d\nu^k,
\nonumber
\\
&&\mathcal{D}_\infty^{\nu^k}=\{f\in {\cal D}_0^k :  \cE_{\nu^k}^{a,k} (f,f)+
\|f \|_{L^2(\R^k\times\mM,\, \nu^k)}^2 <\infty \}.
\nonumber
\end{eqnarray}
It was proved that $(\cE_{\nu^k}^{a,k},{\cal D}_\infty^{\nu^k})$ is closable and its closure $(\cE_{\nu^k}^{a,k},\mathcal{D}^{\nu^k})$ is a quasi-regular Dirichlet form in case $\mu$ is a quasi Gibbs measure \cite{o.tp}. 
Let $\nu_\star^k$ be the Campbell measure associated with DPP $\mu_\star$, $\star\in\{\sin, \Ai\} \cup \{J_\nu\ ;\nu>-1\}$.
For the case where $(\nu^k, a(x))= (\nu_{\sin}^k, 1)$, $ (\nu_{\Ai}^k,1)$, 
and  $(\nu_{J_\nu}^k, 4x)$,
we denote the quasi-Dirichlet form by $(\cE_\star^k, \mathcal{D}^k)$ and the associated diffusion process by
$((\X^k(t),H(t), \P_\star^{\x_k,\eta})$, $\star\in\{\sin, \Ai\} \cup \{J_\nu\ ;\nu>-1\}$.
We introduce the map $\u: \R^k \times \mM \to \mM$ defined as
$$
\u(\x_k,\xi)=\sum_{j=1}^k \delta_{x_j} + \xi, \quad \x_k\in\R^k, \ \xi\in\mM.
$$
It is shown in \cite{o.tp} that if $\u(\x_k, \eta)=\xi$, process $(\u(\X^k(t),H(t)), \mathbf{P}_\star^{\x_k,\eta})$ coincides with process $(\Xi(t), \mathbf{P}_\star^{\xi}$). 

Let $\star \in \{\sin, \Ai\} \cup \{J_\nu\ ;\nu>-1\}$.
Theorem \ref{TH:SM} implies that the Dirichlet form $(\hat{\cE}_{\star}, \hat{\cD}_{\star})$ associated with process $(\Xi(t), \P_{\star})$  is quasi-regular \cite{mr, FOT}.
A function $f$ on the configuration space $\mM$ is said to be
polynomial if it is written in the form
\begin{equation}\label{pol}
f(\xi)=F \left(
\int_{\R}\phi_1(x)\xi(dx), \int_{\R}\phi_2(x)\xi(dx),
\dots,\int_{\R}\phi_k(x)\xi(dx) \right)
\end{equation}
with polynomial function $F$
on ${\R}^{k}, k \in \N$,
and smooth functions $\phi_j, 1 \leq j \leq k$
on $\R$ with compact supports.
Let $\cP$ be the set of all polynomial functions on $\mM$.
In \cite[Proposition 7.2]{KT07a} we showed that 
\begin{equation}
\hat{\cE}_\star(f,g) = \cE_\star(f,g),
\quad f,g \in \cP.
\label{eq:Pol}
\end{equation}
$(\cE_\star, \cD_\star)$ and $(\hat{\cE}_\star, \hat{\cD}_\star)$ are then closed extensions of $(\cE_\star, \cP)$, and the former is the smallest one \cite{ot.0}.
These relations are generalized to $k$-labeled dynamics.

\begin{lem} \label{L:38}
Let $\star\in\{\sin, \Ai\} \cup \{J_\nu\ ;\nu>-1\}$. For each $k\in \N$, there exists a diffusion process $((\X^k(t), H(t)), \hat{\P}_{\star}^{\x_k,\eta})$
associated with quasi-regular Dirichlet form $(\hat{\cE}_{\star}^{k}, \hat{\cD}_{\star}^k)$ such that
\begin{eqnarray}
&&(\u(\X^k(t),H(t)), \mathbf{P}_{\star}^{\x_k,\eta})=(\Xi (t),\P_\star^{\xi}), \quad \mbox{if  $\u(\x_k,\eta)=\xi$},
\label{37:1}
\\
&&\cE_{\star}^k(f,g)=\hat{\cE}_{\star}^k(f,g)
\quad f,g \in C_0(\R^k)\otimes {\cal P}.
\label{37:2}
\end{eqnarray}
\end{lem}

\noindent Proof. Let $\star\in \{\sin, \Ai\} \cup \{J_\nu\ ;\nu>-1\}$.
We introduce $(\tilde{\cE}_{\star}^{k}, \tilde{\cD}_{\star,0}^k)$ defined as
$$
\tilde{\cD}_{\star,0}^k=\{ f(\x_k,\eta)= g(\x_k,\u(\x_k,\eta)) ; g \in C_0^{\infty}(\R^k)\otimes \cD_{\star} \},
$$
$$
\tilde{\cE}_{\star}^k (f,g)=\int_{\mM}\nabla^k[f,g](\x_k,\eta)+\D_\star^k [f(\x_k,\cdot),g(\x_k,\cdot)](\eta) \nu_{\star}^k(d\x_k d\eta),
\quad f,g \in \tilde{\cD}_{\star,0}^k,
$$
where the derivatives are taken in the sense of the Schwartz distribution, $\D_{\star}^k=\D^{1,k}$ if $\star\in\{\sin,\Ai\}$, and $\D_{\star}^k=\D^{4x,k}$ if $\star\in\{J_\nu, \nu>-1\}$.
For $\eta\in \mM$ and $r\in\N$, we set $\eta_r=\eta_{[-r,r]}=\sum_{j=1}^\ell \delta_{y_j}$ and $\zeta=\eta_{[-r,r]^c}$.
For $  f \in \tilde{\cD}_{\star,0}^k $ we set
$$
{f}_{r,\zeta}(\x_k,\eta_r)= f(\x_k,\eta)\1_{[-r,r]^k}(\x_k)
$$
and 
\begin{eqnarray}
&&\D_{\star,r}^{k}[f,g][\x_k,\eta]=\nabla^k[f_{r,\zeta},g_{r,\zeta}](\x_k,\eta_r)+\D_{\star}^k [f_{r,\zeta}(\x_k,\cdot),g_{r,\zeta}(\x_k,\cdot)](\eta_r)
\nonumber\\
&&\tilde{\cE}_{\star, r}^{k} (f,g)= \int_{\mM}{\Bbb D}_{\star,r}^{k}(f,g)d\nu^k, \quad f,g\in \tilde{\cD}_{\star,0}^k.
\nonumber
\end{eqnarray}
It is readily seen that the bilinear forms $(\tilde{\cE}_{\star, r}^{k}, \tilde{\cD}_{\star,0}^k)$, $r\in\N$ are closable and increasing. From \cite[Lemma 2.1 (1)]{Osa96} we see that $(\tilde{\cE}_{\star}^{k}, \tilde{\cD}_{\star,0}^k)$ is closable. 
The quasi-regularity of the closure $(\hat{\cE}_{\star}^{k}, \hat{\cD}_{\star}^k)$  of $(\tilde{\cE}_{\star}^{k}, \tilde{\cD}_{\star}^k)$ is derived from that of the Dirichlet form $(\hat{\cE}_\star, \hat{\cD}_\star )$. Hence, the associated diffusion process $((\X^k(t), H(t)), \P_{\star}^{\x_k,\eta})$
can be constructed. 
Equation (\ref{37:2}) is derived from (\ref{eq:Pol}), while (\ref{37:1}) is derived from the argument to show \cite[Lemma 4.2]{o.tp}. This completes the proof.
\qed

\vskip 3mm

\noindent {\it Proof of Theorem \ref{TH:SDE}}.
We show Theorem \ref{TH:SDE} (i) by applying \cite[Theorem 26]{o.isde} to process $(\Xi(t), \P_{\star})$ associated with the Dirichlet form $(\hat{\cE}_{\star}, \hat{\cD}_{\star})$, $\star\in\{\sin, \Ai\} \cup \{J_\nu\ ;\nu>-1\}$. Of the assumptions (A.1)--(A.5) of the theorem,  (A.1), (A.2), and (A.5) are satisfied, because the related measures $\mu_\star$ are the same as those of the Dirichlet form $(\cE_{\star}, \cD_{\star})$, which have already been verified in \cite{o.isde, ho.1, ot.1}.
Assumption (A.4) is derived from the fact that the capacity related to $(\cE_{\star}, \cD_{\star})$ is greater than that of $(\hat{\cE}_{\star}, \hat{\cD}_{\star})$, since these Dirichlet forms are closed extensions of $(\cE_\star, \cP)$, and the former is the smallest one.

Condition (A.3) is used in the proof of \cite[Theorem 26]{o.isde} to construct a $k$-labeled process and to check that process $X_j(t)-X_j(0)$, $j=1,2,\dots, k$ can be regarded as a Dirichlet process. These claims are derived from Lemma \ref{L:38} together with the fact that $\cD_\star^k \subset \hat{\cD}_\star^k$.
(i) is thus proved. 

Claim (ii) is derived from (i) and
the uniqueness of the strong solutions of ISDEs (\ref{ISDE:sin}), (\ref{ISDE:Ai}), and (\ref{ISDE:J}), which have been proved in \cite{ot.2}.
\qed

\vskip 3mm

\noindent {\it Proof of Corollary \ref{C:23}.}
Process $(\Xi(t), \mathbf{P}_{\star}^{\xi^N})$, $\star\in\{\sin, \Ai\} \cup \{J_\nu\ ;\nu>-1\}$ is identical in distribution to process $(\Xi(t), \tilde{\P}_{\star}^{\xi^N})$. Corollary \ref{C:23} (i) is thus derived from Theorem \ref{TH:SDE} and
Proposition \ref{Prop:34}. Claim (ii) is readily derived from Claim (i).\qed
\vskip 10mm

\begin{small}

\noindent{\bf Acknowledgements} \
H.O. is supported in part by a Grant-in-Aid for Scientific Research (KIBAN-A, No. 24244010)
of the Japan Society for the Promotion of Science.
H.T. is supported in part by a Grant-in-Aid for Scientific Research (KIBAN-C, No. 23540122)
of the Japan Society for the Promotion of Science.


\end{small}
\end{document}